\documentclass[graybox]{svmult}

\usepackage{mathptmx}       
\usepackage{helvet}         
\usepackage{courier}        

\usepackage{makeidx}         
\usepackage{graphicx}        
\usepackage{multicol}        
\usepackage[bottom]{footmisc}

\usepackage{cite}
\usepackage{fancyvrb}        

\usepackage{amsmath}
\usepackage{tikz}
\usetikzlibrary{arrows.meta, spy}

\begin{document}

\title*{Isogeometric Simulation and Shape Optimization with Applications to Electrical Machines}
\titlerunning{Isogeometric Simulation and Shape Optimization for Electrical Machines}
\titlerunning{Isogeometric Simulation and Shape Optimization}
\author{Peter Gangl \inst{1}
\and
Ulrich Langer \inst{2,4}
\and
Angelos Mantzaflaris \inst{3}
\and
Rainer Schneckenleitner \inst{4}
}
\institute{
	Peter Gangl \at Institute of Applied Mathematics,  TU Graz,  Steyrergasse 30, 8010 Graz, Austria
\email{gangl@math.tugraz.at}
  \and Ulrich Langer \at Institute of Computational Mathematics, JKU Linz, Altenberger Straße 69, 4040 Linz, Austria
  \email{ulanger@numa.uni-linz.ac.at}
  \and Angelos Mantzaflaris \at Institute of Applied Geometry, JKU Linz, Altenberger Straße 69, 4040 Linz, Austria
  \email{angelos.mantzaflaris@jku.at}
  \and Rainer Schneckenleitner \at RICAM, Austrian Academy of Sciences, Altenberger Straße 69, 4040 Linz, Austria
  \email{rainer.schneckenleitner@ricam.oeaw.ac.at}
  }

%
%
\maketitle

\abstract{
	Future e-mobility calls for efficient electrical machines.
	For different areas of operation, these machines have to satisfy certain desired properties
	that often depend on their design.
	Here we investigate the use of multipatch Isogeometric Analysis (IgA) for the simulation and shape optimization  of the electrical machines.
	In order to get fast simulation and optimization results,
	we use non-overlapping domain decomposition (DD) methods to solve the large systems of algebraic
	equations arising from the IgA discretization of underlying partial differential equations.
	The DD is naturally related to the multipatch representation of the computational domain, and provides the framework for the parallelization of the DD solvers.
	}

\section{Introduction}
\label{Schneckenleitner:sec_introduction}

Isogeometric Analysis (IgA) is a relatively new approach
which was introduced in \cite{GLMS:Hughes}.
It can be seen as an alternative to the more classical Finite Element Method (FEM).
The idea in IgA is to use the same basis functions for both representing the computational geometry and  solving the  partial differential equations (PDEs). This aspect makes IgA especially interesting for design optimization procedures.
In practice, it is often the case that one
performs
design optimization and geometric modeling simultaneously.
State-of-the-art computer aided design (CAD) software
use
B-splines or NURBS for the modeling process
whereas the design optimization requires an analysis suitable representation of the model.
So far the design optimization is mainly done using
FEM as discretization method.
Hence, the B-spline or NURBS representation of the geometric model has to be converted into
a
suitable mesh for the
Finite Element Analysis.
This conversion is in general very computationally demanding.
The new IgA paradigm circumvents these problems. Therefore, IgA is very beneficial for the simulation and optimization
when the representation of the computational domain comes from CAD software;
see \cite{GLMS:Schneckenleitner,GLMS:Bontnick} for applications to electrical machines.

Since practical optimization problems tend to be very large,
the numerical solution of the underlying PDEs becomes computationally very expensive.
Moreover, in PDE-constrained shape optimization processes, there are more than one
PDE to solve.
In particular, line search requires to solve the magnetostatic PDE constraint several times.
In order to get fast optimization results, we use
Dual-Primal IsogEometric Tearing and Interconnecting (IETI-DP) methods
for the solution of the linear algebraic systems
arising from the IgA discretization.
The IETI-DP solvers are  non-overlapping domain decomposition methods;
see \cite{GLMS:KleissPechsteinJuettlerTomar:2012a,GLMS:Hofer}.
IETI-DP methods are closely related the FEM-based FETI-DP methods; see, e.g., \cite{GLMS:Pechstein:2013a} and the references therein.
We show
that  IETI-DP methods are superior to sparse direct solvers with respect to computational time
and memory requirement.
Moreover, IETI-DP provides a natural framework for parallelization.
Indeed, our numerical experiments on a distributed memory computer show an excellent scaling
behavior of this method.


\section{Shape optimization via gradient descent}
\subsection{Problem Description}
\label{Schneckenleitner:Simulation_and_Optimization}
\begin{figure}[ht]
 	\centering
 	\includegraphics[width = 5.7cm, height = 5.7cm]{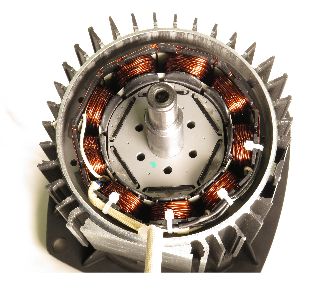}
 	\includegraphics[width = 5.2cm, height = 5.2cm]{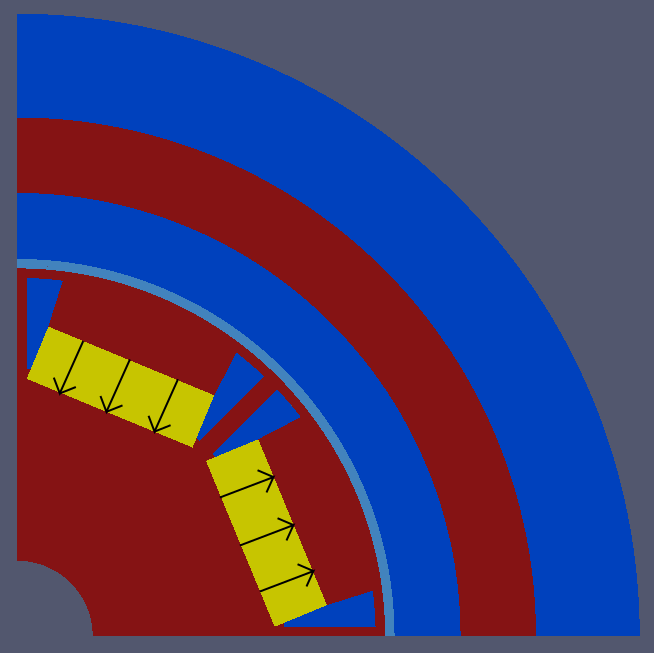}
 	\caption{Real world IPM electric motor and a model of a quarter of its cross section.}
 	\label{fig::IPM+ComputationalModel}
 \end{figure}

We investigate the simulation and shape optimization of an interior permanent magnet (IPM)
electric motor by means of IgA. The IgA approach seems to be very attractive for such practical problems.
The most beneficial aspect of IgA in the context of optimization is the fact that the same basis functions which are used to represent the geometry of the IPM electric motor are also exploited to solve the underlying PDEs. In the optimization procedure, we want to optimize the shape of the motor in order to maximize the runout performance, i.e. to maximize the smoothness of the rotation of the motor. An example of an IPM electric motor is given in Fig.~\ref{fig::IPM+ComputationalModel} (left). One possible way to optimize the runout performance of an IPM electric motor is to minimize the squared $L^2$-distance between the radial component of the magnetic flux in the air gap and a desired smooth reference function $B_\mathrm{d}$. The resulting optimization problem is subject to the 2d magnetostatic PDE as constraint.

 Mathematically, the arising optimization problem can be expressed as follows:
 \begin{eqnarray}
 \label{Schneckenleitner:objective}
 \underset{D}{\mbox{min}}\; J(u):= \int_{\varGamma}^{} |B(u) \cdot \mathrm{n}_\varGamma - B_\mathrm{d}|^2 \mathrm{d}s = \int_{\varGamma}^{} |\nabla u \cdot \mathrm{\tau}_{\varGamma} - B_\mathrm{d}|^2 \mathrm{d}s \\
 \mbox{s.t.} \; u \in H_0^1(\Omega): \int_{\Omega}^{} \nu_\mathrm{D}(x)
 \nabla u \cdot \nabla \eta \; \mathrm{d}x
 = \langle F, \eta \rangle \; \forall \eta \in H_0^1(\Omega),
 \label{Schneckenleitner:constraint}
 \end{eqnarray}
where $J$ denotes the objective function, $\varGamma$ is the midline of the air gap, $\Omega$ denotes the whole computational domain,
 and $D$ is the domain of interest also called design domain. The variational problem (\ref{Schneckenleitner:constraint})
 is nothing but
 the 2d linear magnetostatic problem with the piecewise constant magnetic reluctivity $\nu_D(x)=\chi_{\Omega_f(D)}(x)\nu_1 + \chi_{\Omega_{\text{mag}}}(x) \nu_{\text{mag}} + \chi_{\Omega_{\text{air}}(D)}(x)\nu_0$. Here, $\Omega_f$, $\Omega_{\text{mag}}$ and $\Omega_{\text{air}}$ denote the ferromagnetic, permanent magnet and air subdomains, respectively, and $\nu_1$, $\nu_{\text{mag}}$ and $\nu_0$ denote the corresponding reluctivity values. Note that the shape $D$ enters the optimization problem via the function $\nu_D$ and influences the objective function via the solution $u$. The right hand side $F \in H^{-1}(\Omega)$ in (\ref{Schneckenleitner:constraint}) is defined by the linear functional
 \begin{equation}
	\langle F, \eta\rangle := \int_{\Omega}^{} (J_3\eta + \nu_\mathrm{mag} M^\bot \cdot \nabla \eta) \; \mathrm{d}x
 \end{equation}
 for all $\eta \in H_0^1(\Omega)$. Here, $M^\bot$ denotes the perpendicular of the magnetization $M$, which is indicated in Fig.~\ref{fig::IPM+ComputationalModel} and vanishes outside the permanent magnets, and $J_3$ is the third component of the impressed current density in the coils. Note that the solution $u$ is the third component of the magnetic vector potential, i.e. $B(u) = \mbox{curl}((0,0,u)^T)$. Moreover, $n_\Gamma = (n_1, n_2, 0)^T$ and $\tau_\Gamma = (\tau_1, \tau_2)^T$ denote the outward unit normal and unit tangential vectors along the air gap, respectively.

 We are interested in the radial component of the magnetic flux density along the air gap due to the permanent magnetization.
 For that reason, we set $J_3=0$ and consider the coil regions as air. Figure~\ref{fig::IPM+ComputationalModel} (right) shows a quarter of a cross section of a simplified IPM electric motor that is provided by CAD software.
 Hence, this geometry representation is suitable for IgA simulation.
 The red-brown areas represent ferromagnetic material ($\Omega_f$), the blue areas consist of air ($\Omega_{\text{air}}$), the yellow areas are the permanent magnets ($\Omega_{\text{mag}}$).
 The air gap of the motor is highlighted in light blue.
 In this initial model for the optimization, the design domain $D$ is the ferromagnetic area right above the permanent magnets. In order to get a smoother rotation we are looking for a better shape of this part $D$.

%
%

\subsection{The shape derivative}
\label{Schneckenleitner:shape_calculus}
For the optimization of the IPM electric motor, we use gradient based optimization techniques.
Hence, we need the derivative of the objective $J$ with respect to a change of the current shape. The shape derivative in tensor form
\cite{GLMS:Gangl, GLMS:Schneckenleitner, GLMS:Sturm} of our optimization problem is given by
\begin{equation}
\label{Schneckenleitner:derivative}
dJ(D)(\phi) = \int_{\Omega}^{}\mathcal{S}(D, u, p) : \partial \phi \mathrm{d}x,
\quad \forall \phi \in H^1_0(\Omega, \rm{I\!R}^2)
\end{equation}
with
$
\mathcal{S}(D, u, p) = (\nu_\mathrm{D}(x)\nabla u \cdot \nabla p - \nu_{\text{mag}} \nabla p \cdot M^\bot) \mathcal{I}
                        + \nu_{\text{mag}} \nabla p \otimes M^\bot
                        - \nu_\mathrm{D}(x) \nabla p \otimes \nabla u
                        - \nu_\mathrm{D}(x) \nabla u \otimes \nabla p,
$
where $\mathcal{I}$ denotes the identity, the state u solves the constraint $(\ref{Schneckenleitner:constraint})$, and $p$ solves the adjoint problem
\begin{equation}
\label{Schneckenleitner:adjoint}
\int_{\Omega}^{} \nu_\mathrm{D}(x)
\nabla p \cdot \nabla \eta \; \mathrm{d}x
=
-2\int_{\varGamma}^{} (B(u) \cdot n_{\varGamma} - B_\mathrm{d})(B(\eta) \cdot n_{\varGamma}) \; \mathrm{d}s
\quad \forall \eta \in H^1_0(\Omega).
\end{equation}

%
%

\subsection{Numerical shape optimization}

\begin{figure}[ht]
	\centering
	\includegraphics[width = 5.2cm, height = 5.2cm]{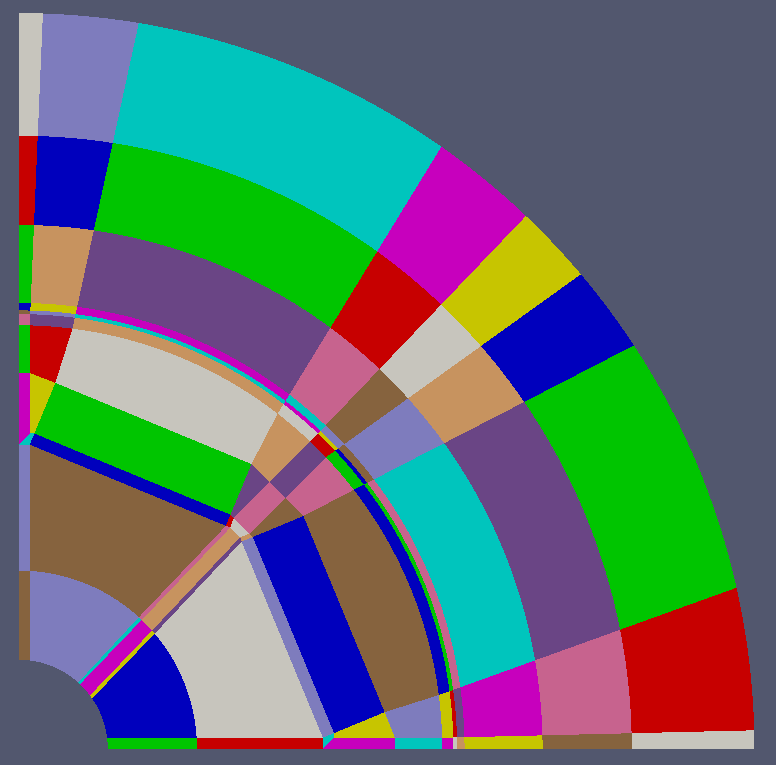}\hspace*{5mm}
	\includegraphics[width = 5.2cm, height = 5.2cm]{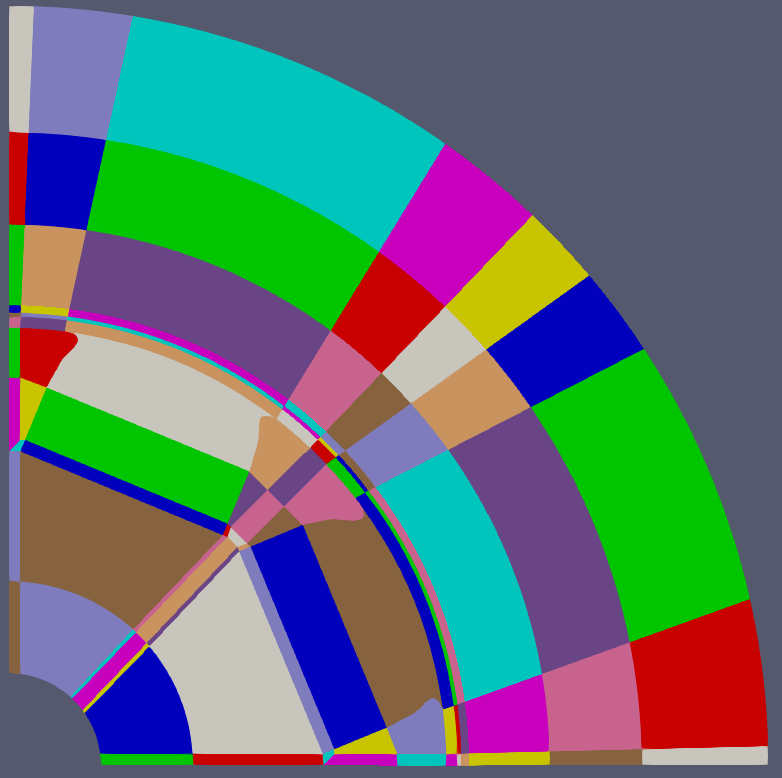}
	\caption{Initial and final design of an IPM motor.}
	\label{Schneckenleitner:Initial+FinalIPM}
\end{figure}

We  used a continuous Galerkin (cG) IgA discretization for both the simulation and optimization problems.
The implementation is done in \textbf{G+Smo}\footnote{Mantzaflaris, A. et al.: G+Smo (geometry plus simulation modules) v0.8.1., \url{http://gs.jku.at/gismo, 2017  Jun 19 2018}}.
Figure \ref{Schneckenleitner:Initial+FinalIPM} (left) shows a possible computational domain suitable for cG.
The shown multipatch domain consists of 93 patches. For each of these patches, we used a B-spline mapping from a reference patch with splines of degree 3.
For the optimization, we need the shape gradient $\nabla J \in V:= H^1_0(\Omega, {\rm I\!R}^2)$ which can be computed by solving the auxiliary problem: find $\nabla J \in V$ such that \begin{equation}
   \label{Schneckenleitner:auxiliary}
   b(\nabla J, \psi) = -\mathrm{d}J(D)(\psi) \quad \forall \psi \in V.
\end{equation}
The expression on the right hand side of (\ref{Schneckenleitner:auxiliary}) is the negative shape derivative whereas the expression $b(\cdot, \cdot)$ on the left hand side is some $V$-elliptic, $V$-bounded bilinear form which must be chosen appropriately. For our studies, we used
\begin{equation}
   \label{Schneckenleitner:auxiliaryBiform}
   b(\phi, \psi) = \int_{\Omega}^{} \phi \cdot \psi \; \mathrm{d}x + \int_{\Omega}^{} \alpha (\partial \phi : \partial \psi) \; \mathrm{d}x
\end{equation}
with a patchwise constant function $\alpha \in L^\infty(\Omega)$.

In the right picture of Fig.~\ref{Schneckenleitner:Initial+FinalIPM}, we can see the optimized shape with respect to the runout performance compared to the initial domain on the left. We were able to reduce the objective from
$4.236 \cdot 10^{-4}$ down to $2.781 \cdot 10^{-4}$.


\section{Fast numerical solutions by IETI-DP}

\begin{figure}[ht]
	\centering
	\includegraphics[width = 5.2cm, height = 5.2cm]{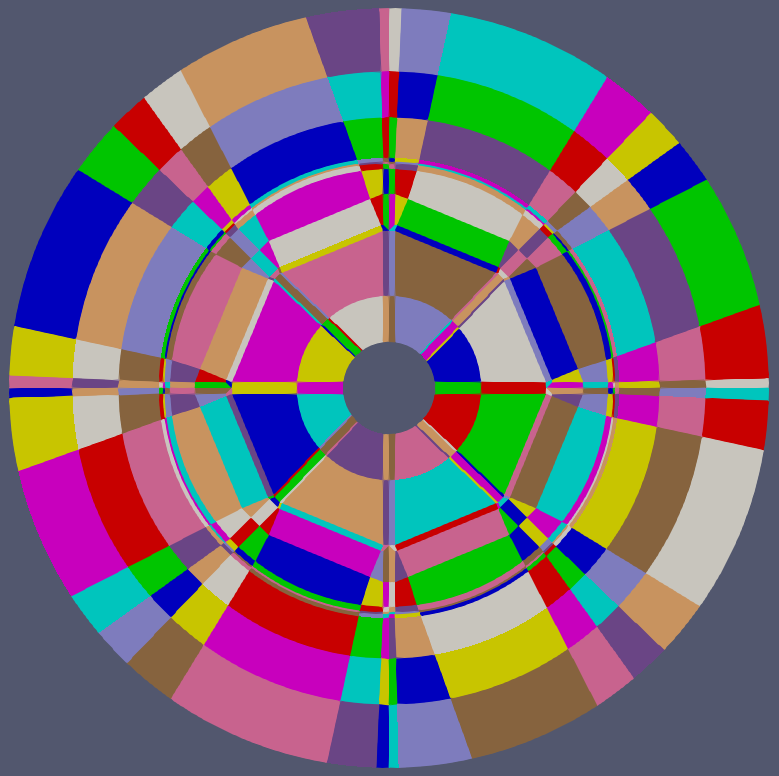}\hspace*{5mm}
	\includegraphics[width = 5.2cm, height = 5.2cm]{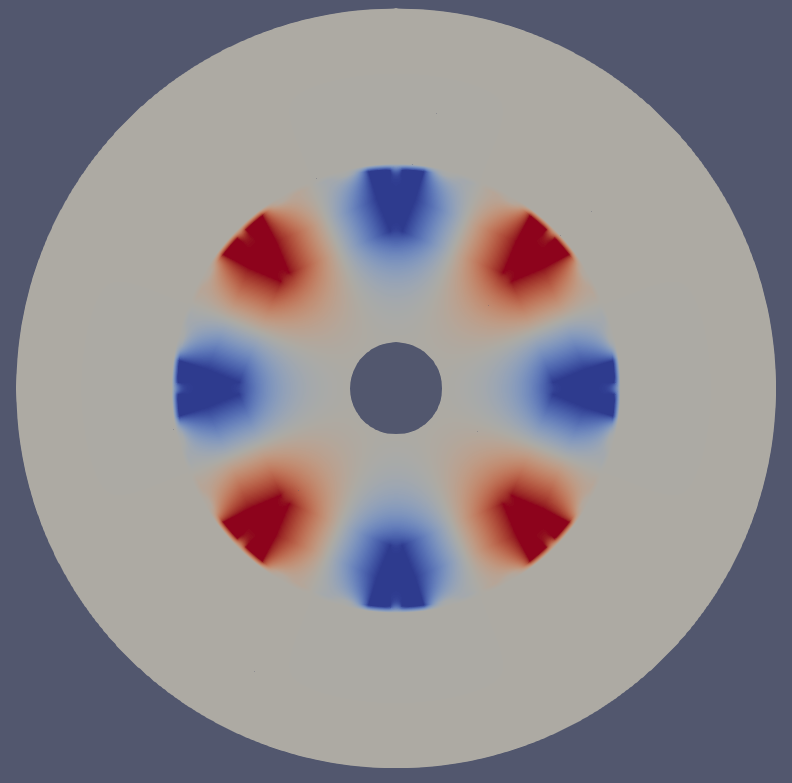}
	\caption{Whole initial cross section as well as the solution.}
	\label{Schneckenleitner:Initial+solState}
\end{figure}

Up to now, we have solved the arising PDEs
by means of a sparse direct solver.
One drawback of  a direct solution method is that it is rather slow
for large-scale systems.
In particular,
in shape optimization,
we have to solve the state equation (\ref{Schneckenleitner:constraint}), the adjoint equation (\ref{Schneckenleitner:adjoint}), and the auxiliary problem (\ref{Schneckenleitner:auxiliary}) for the shape gradient, which decouples into two scalar problems, in every iteration of the optimization algorithm. Moreover, during a line search procedure, it might be the case that the state equation has to be solved several times. To overcome the issue of a slow performance, we were looking for a fast and suitable solver for our simulation and optimization processes. We chose the IETI-DP technique for solving the PDEs.
IETI-DP is a  non-overlapping domain decomposition technique which introduces local subspaces which are then again coupled using additional constraints. A comparison between a sparse direct solver and IETI-DP for solving the state equation (\ref{Schneckenleitner:constraint}) on a full cross section of an IPM electric motor clearly shows that the recently developed IETI-DP method \cite{GLMS:Hofer} performs much better as can be seen in
Table~\ref{Schneckenleitner:speedup}. This and the subsequent numerical experiments were done on
RADON1 ({\url{https://www.ricam.oeaw.ac.at/hpc/overview/}})
a high performance computing cluster with 1168 computing cores and 10.7 TB of memory.
Table~\ref{Schneckenleitner:speedup} also shows that,
with an increasing number of degrees of freedom,
the proposed IETI-DP technique solves the problem much faster
than the sparse  direct solver.
Moreover, it can be seen that, with too many degrees of freedom, the sparse direct solver ran out of memory whereas IETI-DP could provide the solution to the problem.
The solution to the state equation is shown in Fig.~\ref{Schneckenleitner:Initial+solState} (right).

\begin{table} [h!]
	\caption{SuperLU  vs. IETI-DP on a single core.}
	\label{Schneckenleitner:speedup}
	\renewcommand{\arraystretch}{1.5}
	\begin{center}
		\begin{tabular}{r|r|r|c}
			\# dofs & \multicolumn{1}{c|}{SuperLU} & \multicolumn{1}{c|}{IETI-DP} & speedup \\
			\hline \hline
			72 572		& 36.0 sec	& 17.0 sec	& 2.12	\\
			250 844		& 193.0 sec	& 69.8 sec	& 2.77	\\
			928 796		& 1943.0 sec	& 463.0 sec	& 4.20  \\
			3 570 332	&  -- \hspace*{6mm}	& 1179.0 sec	& --
		\end{tabular}
	\end{center}
\end{table}

Moreover, IETI-DP provides a natural framework for parallelization.
Because of the multipatch structure of the computational domains in IgA,
each patch can be seen as a subdomain in the IETI-DP approach. Then one can create suitable subdomains consisting of a certain number of patches for each processor, e.g., one possible choice is to group the patches to subdomains according to their number of degrees of freedom which means that the degrees of freedom are almost evenly distributed over the number of processors.
Table \ref{Schneckenleitner:strongScaling} shows the strong scaling behavior of the IETI-DP
solver.
In this experiment, we solved the constraint equation (\ref{Schneckenleitner:constraint}) on the full cross section of an IPM electric motor with 3 570 332 degrees of freedom.
From Table \ref{Schneckenleitner:strongScaling}, we can see the expected performance,
i.e., if we double the number of processors the computation time reduces nearly by a factor of two.

\begin{table} [h]
	\caption{Strong scaling with IETI-DP and 3 570 332 dofs.}
	\label{Schneckenleitner:strongScaling}
	\renewcommand{\arraystretch}{1.5}
	\begin{center}
		\begin{tabular}{c||c|c|c|c|c|c|c|c}
			\# cores & 1 & 2 & 4 & 8 & 16 & 32 & 64 & 128 \\
			\hline
			time [sec]		& 1179 	& 577 	& 325 & 164 & 89 & 43 & 22 & 14  \\
			\hline
			rate     	& -- & 2.04 & 1.78 & 1.98 & 1.84 & 2.07 &  1.95 & 1.57
		\end{tabular}
	\end{center}
\end{table}


\section{Shape optimization based on Ipopt and IETI-DP}
In this section we point out the usage of \textbf{Ipopt},
which stands for \textbf{I}nterior \textbf{p}oint \textbf{opt}imizer \cite{GLMS:Ipopt},
for the shape optimization using IETI-DP as
underlying PDE solver.
If we do shape optimization without any additional considerations,
then we might run into troubles.
More precisely,
it can happen that we get self-intersections in the final shape even if the objective decreases.

To prevent such self-intersections,
we consider the Jacobian determinant of the geometry transformation in the design domain and its neighboring air regions. The Jacobian determinant of these patches must have the same sign in each iteration. If the sign changes from one iteration to the next, then we reduce the step size until the Jacobian determinant of the new design has the same sign as in the initial configuration.
In this way, we are able to ensure that the shape is technically feasible.

In the first naive approach, all control coefficients of the multipatch domain are considered as design variables, and the vector field computed by (\ref{Schneckenleitner:auxiliary}) is applied globally. The computational effort for the optimization can be reduced by applying the computed vector field only on the important interfaces between the design domain and the neighboring air regions. This reduces the number of design variables from approximately $28 000$ to $128$ in the coarsest setting. The inner control coefficients of the design area and the bordering air regions are rearranged via a spring patch model \cite{GLMS:Manh}.

In a first test setting, Ipopt stops at an optimal solution after 95 iterations
using a BFGS method.
We set the NLP error tolerance to $10^{-6}$, the relative error in the objective change to the same value,
and we decided to exit the optimization loop after three iterations within these error bounds.
The objective value dropped from $4.266 \cdot 10^{-4}$ down to  $2.587 \cdot 10^{-4}$

\begin{figure}[h]
	\centering
	\begin{tikzpicture}
		\begin{scope}[spy using outlines={draw,yellow,height=2.5cm, width=4.2cm,magnification=4, connect spies}]
		\node[anchor=south west,inner sep=0] at (0,0) {\includegraphics[interpolate = true, width=5.2cm, height = 5.2cm]{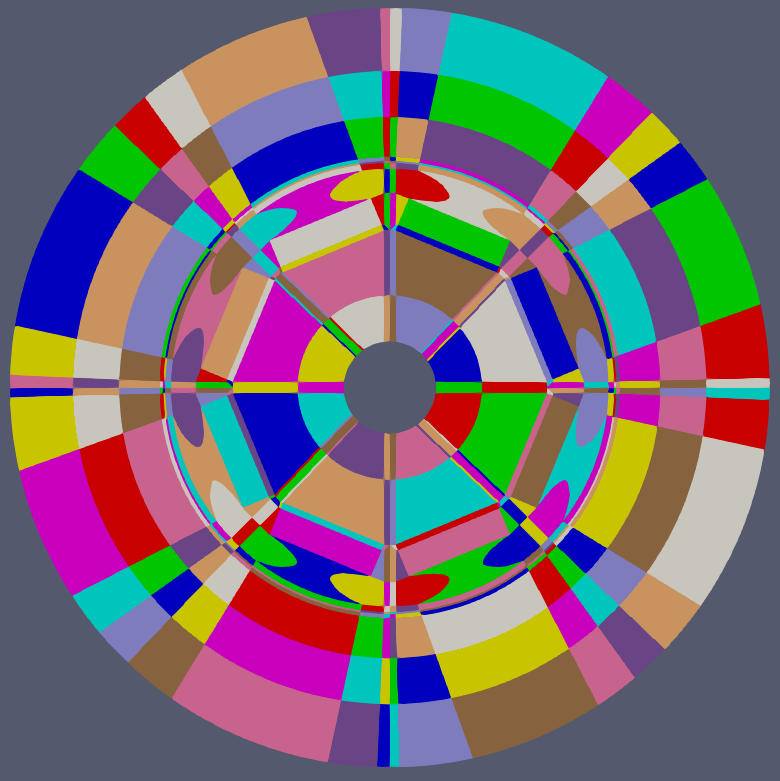}};

		\spy on (3.12,3.85) in node [left] at (10,3.25);
		\end{scope}
    \end{tikzpicture}
	\caption{
	Optimal shape after $130$ iterations with relaxed bounds (left), zoom into one of the design regions (right)}
	\label{Schneckenleitner:Optimal130+DesDomReparam130}
\end{figure}

Furthermore, we tried an additional experiment were we relaxed the bounds on the constraints a bit.
In particular, we set the \texttt{bound{\_}relax{\_}factor} in Ipopt to $1$.
The result of this experiment can be seen in Fig.~\ref{Schneckenleitner:Optimal130+DesDomReparam130}.
We may observe from Fig.~\ref{Schneckenleitner:Optimal130+DesDomReparam130}
that we get a very smooth final shape with even a smaller objective value of $2.436 \cdot 10^{-4}$.
We point out that if we adjust the different optimization parameters we may get different optimal shapes and different objective values in the end.

\begin{acknowledgement}
This work was
supported by the Austrian Science Fund  (FWF)
via the grants NFN S117-03 and the DK W1214-04.
We also acknowledge the permission to use the Photo in Fig.~\ref{fig::IPM+ComputationalModel} (left) taken by the Linz Center of Mechatronics (LCM). The  motor was produced by Hanning Elektro-Werke GmbH \& Co KG.
\end{acknowledgement}


\end{document}